\documentclass[aos,MSNbibl,nameyear,dvips]{arximspdf}

\doi{10.1214/13-AOS1084} %kopijuoti is PTS
\volume{41}
\issue{1}
\pubyear{2013}
\firstpage{177}
\lastpage{195}

\makeatletter
\newcommand{\rrvert}{\vert}
\newcommand{\llvert}{\vert}
\def\cal{\mathcal}
\newproclaim{example}{Example}
\newtheorem{theorem}{Theorem}
\makeatother

\begin{document}
\begin{frontmatter}

\title{The subset argument and consistency of MLE in GLMM:
Answer to an open problem and beyond}
\runtitle{Consistency of MLE in GLMM}

\begin{aug}
\author{\fnms{Jiming} \snm{Jiang}\corref{}\ead[label=e1]{jiang@wald.ucdavis.edu}\thanksref{t1}}
\thankstext{t1}{Supported in part by NIH Grant R01-GM085205A1 and
NSF Grants SES-9978101,
DMS-02-03676, DMS-04-02824, DMS-08-06127 and SES-1121794.}
\runauthor{J. Jiang}
\affiliation{University of California, Davis}
\address{Department of Statistics\\
University of California, Davis\\
Davis, California 95656\\
USA\\
\printead{e1}}
\end{aug}

% HISTORY:
\received{\smonth{3} \syear{2012}}
\revised{\smonth{11} \syear{2012}}

% ABSTRACT
%
\begin{abstract}
We give answer to an open problem regarding consistency of the
maximum likelihood estimators (MLEs) in generalized linear mixed models
(GLMMs) involving crossed random effects. The solution to the open problem
introduces an interesting, nonstandard approach to proving consistency of
the MLEs in cases of dependent observations. Using the new technique, we
extend the results to MLEs under a general GLMM. An example is used to
further illustrate the technique.
\end{abstract}

% KEYWORDS
% Pirmas kwd is didziosios raides
%


\end{frontmatter}
%
%s1 #&#
\section{Introduction}\label{sec1}

Generalized linear mixed models (GLMMs) have become a popular and very
useful class of statistical models. See, for example, \citet{Jia07},
McCulloch, Searle and Neuhaus (\citeyear{McCSeaNeu08}) for some wide-ranging accounts of
GLMMs with theory and applications. In the earlier years after GLMM was
introduced, one of the biggest challenges in inference about these
models was computation of the maximum likelihood estimators (MLEs).
As is well known, the likelihood function under a GLMM typically
involves integrals that cannot be computed analytically. The
computational difficulty was highlighted by the infamous salamander
mating data, first introduced by McCullagh and Nelder [(\citeyear{McCNel83}), Section 14.5].
A mixed logistic model, which is a special case of GLMM, was proposed
for the salamander data that involved crossed random effects for the
female and male animals. However, due to the fact that the random effects
are crossed, the likelihood function involves a high-dimensional integral
that not only does not have an analytic expression, but is also difficult
to evaluate numerically [e.g., \citet{Jia07}, Section 4.4.3]. For years, the
salamander data has been a driving force for the computational
developments in GLMM. Virtually every numerical procedure that was
proposed used this data as a ``gold standard'' to evaluate, or
demonstrate, the procedure. See, for example, \citet{KarZeg92},
\citet{BreCla93}, \citet{DruMcC93}, \citet{McC94},
\citet{BreLin95}, \citet{LinBre96}, \citet{Jia98}, \citet{BooHob99},
\citet{JiaZha01}, \citet{SutRao03}, and
\citet{Tor12}.
%s1.1 #&#
\subsection{A theoretical challenge and an open problem}\label{sec1.1}

To illustrate the numerical difficulty as well as a theoretical
challenge, which is the main objective of the current paper, let us
begin with an example.

\begin{example}\label{ex1}
A mixed logistic model was proposed by \citet{BreCla93} for the salamander data, and has since been used [e.g.,
\citet{BreLin95,LinBre96,Jia98}]. Some
alternative models, but only in terms of reparametrizations, have been
considered [e.g., \citet{BooHob99}]. \citet{JiaZha01} noted
that some of these models have ignored the fact that a group of
salamanders were used in both the summer experiment and one of the fall
experiments; in other words, there were replicates for some of the
pairs of female and male animals. Nevertheless, all of these models
are special cases of the following, more general setting. Suppose
that, given the random effects $u_{i}, v_{j}, (i,j)\in S$, where $S$
is a subset of ${\cal I}=\{(i,j)\dvtx1\leq i\leq m, 1\leq j\leq n\}$, binary
responses $y_{ijk}$, $(i,j)\in S$, $k=1,\ldots,c_{ij}$ are conditionally
independent such that, with $p_{ijk}=\mathrm{ P}(y_{ijk}=1|u,v)$, we have
$\operatorname{ logit}(p_{ijk})=x_{ijk}'\beta+u_{i}+v_{j}$, where
$\operatorname{ logit}(p)
=\log\{p/(1-p)\}, p\in(0,1)$, $x_{ijk}$ is an known vector of covariates,
$\beta$ is a unknown vector of parameters, and $u, v$ denote all the random
effects $u_{i}$ and $v_{j}$ that are involved. Here $c_{ij}$ is the number
of replicates for the $(i,j)$ cell. Without loss of generality, assume that
$S$ is a irreducible subset of ${\cal I}$ in that $m, n$ are the smallest
positive integers such that $S\subset{\cal I}$. Furthermore, suppose that
the random effects $u_{i}$'s and $v_{j}$'s are independent with
$u_{i}\sim
N(0,\sigma^{2})$ and $v_{j}\sim N(0,\tau^{2})$, where $\sigma^{2},
\tau^{2}$ are unknown variances. One may think of the random effects
$u_{i}$ and $v_{j}$ as corresponding to the female and male animals, as
in the salamander problem. In fact, for the salamander data, $c_{ij}=2$
for half of the pairs $(i,j)$, and $c_{ij}=1$ for the rest of the pairs.
It can be shown [e.g., \citet{Jia07}, page~126; also see Section~\ref{sec4} in the
sequel] that the log-likelihood function for estimating $\beta,\sigma^{2},
\tau^{2}$ involves an integral of dimension $m+n$, which, in particular,
increases with the sample size, and the integral cannot be further
simplified.
\end{example}

The fact that the random effects are crossed, as in Example~\ref{ex1}, presents
not only a computational challenge but also a theoretical one, that is,
to prove that the MLE is consistent in such a model. In contrast,
the situation is very different if the GLMM has clustered, rather
than crossed, random effects. For example, consider the following.

\begin{example}\label{ex2}
Suppose that, given the random effects $u_{1},
\ldots,u_{m}$, binary responses $y_{ij}, i=1,\ldots,m, j=1,\ldots,n_{i}$
are conditionally independent such that, with $p_{ij}=\mathrm{ P}(y_{ij}
=1|u)$, we have $\operatorname{ logit}(p_{ij})=x_{ij}'\beta+u_{i}$, where
$x_{ij}$ is a vector of known covariates, $\beta$ a vector of unknown
coefficients, and $u=(u_{i})_{1\leq i\leq m}$. Furthermore, suppose
that the $u_{i}$'s are independent with $u_{i}\sim N(0,\sigma^{2})$,
where $\sigma^{2}$ is unknown. It is easy to show that the
log-likelihood function for estimating $\beta,\sigma^{2}$ only
involves one-dimensional integrals. Not only that, a major theoretical
advantage of this case is that the log-likelihood can be expressed as
a sum of independent random variables. In fact, this is a main
characteristic of GLMMs with clustered random effects. Therefore,
limit theorems for sums of independent random variables [e.g., \citet{Jia10}, Chapter
6]
can be utilized to obtain asymptotic properties of the MLE.
\end{example}

Generally speaking, the classical approach to proving consistency of
the MLE [e.g., \citet{LehCas98}, Chapter 6; \citet{Jia10}] relies on
asymptotic theory for sum of random variables, independent or not.
However, one cannot express the log-likelihood in Example~\ref{ex1} as a sum
of random variables with manageable properties. For this reason, it
is very difficult to tackle asymptotic behavior of the MLE in the
salamander problem, or any GLMM with crossed random effects, assuming
that the numbers of random effects in all of the crossed factors
increase. In fact, the problem is difficult to solve even for the
simplest case, as stated in the open problem below.

\begin{quote}
{Open problem} [\textit{e.g., Jiang} (\citeyear{Jia10}), \textit{page 541}]:\vspace*{-2pt}
\textit{Suppose that $x_{ijk}'\beta=\mu$, an unknown parameter, $c_{ij}=1$ for all
$i,j$, $S={\cal I}$, and $\sigma^{2}, \tau^{2}$ are known, say, $\sigma^{2}
=\tau^{2}=1$ in Example~\ref{ex1}. Thus, $\mu$ is the only unknown parameter.
Suppose that $m, n\rightarrow\infty$. Is the MLE of $\mu$ consistent?}
\end{quote}
It was claimed [\citet{Jia10}, pages~541, 550] that even for this seemingly
trivial case, the answer was not known but expected to be anything but
trivial.
%s1.2 #&#
\subsection{Origination of the open problem}\label{sec1.2}

The problem regarding consistency of the MLE in GLMMs with crossed random
effects began to draw attention in early~1997. It remained unsolved over
the past 15 years, and was twice cited as an open problem in the
literature, first in Jiang [(\citeyear{Jia07}), page~173] and later in Jiang [(\citeyear{Jia10}),
page~541]. The latter also provided the following supporting evidence for a
positive answer [\citet{Jia10}, page~550].

Let $k=m\wedge n$. Consider a subset of the data, $y_{ii}, i=1,\ldots,k$.
Note that the subset is a sequence of i.i.d. random variables. It follows,
by the standard arguments, that the MLE of $\mu$ based on the subset,
denoted by $\tilde{\mu}$, is consistent. Let $\hat{\mu}$ denote the MLE
of $\mu$ based on the full data, $y_{ij}, i=1,\ldots,m, j=1,\ldots,n$. The
point is that even the MLE based on a subset of the data, $\tilde{\mu}$,
is consistent; and if one has more data (information), one is expected to
do better. Therefore, $\hat{\mu}$ has to be consistent as well.
%s1.3 #&#
\subsection{The rest of the paper}\label{sec1.3}

In Section~\ref{sec2}, we give a positive answer to the open problem as well as the
proof. Surprisingly, the proof is fairly short, thanks to a new,
nonstandard technique that we introduce, known as the \textit{subset
argument}. Using this argument, we are able to establish both Cram\'{e}r
(\citeyear{Cra46}) and \citet{Wal49} types of consistency results for the MLE. It is
fascinating that a 15-year-old problem can be solved in such
a simple way. The new technique may be useful well beyond solving the
open problem---for proving consistency of the MLE in cases of dependent
observations. We consider some applications of the subset argument in
Section~\ref{sec3} regarding consistency of the MLE in a general GLMM. An example
is used in Section~\ref{sec4} to further illustrate the new technique. Remark and
discussion on a number of theoretical and practical issues are offered
in Section~\ref{sec5}.
%s2 #&#
\section{Answer to open problem}\label{sec2}

Throughout this section, we focus on the open problem stated in Section
\ref{sec1}. Let $\mu$ denote the true parameter.

\begin{theorem}[(Cram\'{e}r consistency)]\label{th1}
There is, with probability tending
to one, a root to the likelihood equation, $\hat{\mu}$, such that
$\hat{\mu}\stackrel{\mathrm{ P}}{\longrightarrow}\mu$.
\end{theorem}

\begin{pf}
The idea was actually hinted in Jiang [(\citeyear{Jia10}), page~550] as
``evidence'' that supports a positive answer (see the last paragraph of
Section~\ref{sec1.2} of the current paper). Basically, the idea suggests that,
perhaps, one could use the fact that the MLE based on the subset data is
consistent to argue that the MLE based on the full data is also
consistent. The question is how to execute the idea. Recall that, in the
original proof of Wald [(\citeyear{Wal49}); also see \citet{Wol49}], the focus was on
the likelihood ratio $p_{\theta}(y)/p_{\theta_{0}}(y)$, and showing that
the ratio converges to zero outside any (small) neighborhood of
$\theta_{0}$, the true parameter vector. Can we execute the subset idea
in terms of the likelihood ratio? This leads to consideration of the
relationship between the likelihood ratio under the full data and that
under the subset data. It is in this context that the following
\textit{subset inequality} (\ref{eq2}) is derived (see Section~\ref{sec5.1} for further
discussion), which is the key to the proof.

Let $y_{[1]}$ denote the (row) vector of $y_{ii}, i=1,\ldots,m\wedge n$,
and $y_{[2]}$ the (row) vector of the rest of the $y_{ij}, i=1,\ldots,m,
j=1,\ldots,n$. Let $p_{\mu}(y_{[1]},y_{[2]})$ denote the probability mass
function (p.m.f.) of $(y_{[1]},y_{[2]})$, $p_{\mu}(y_{[1]})$ the p.m.f. of
$y_{[1]}$,
%
%e1 #&#
\begin{equation}\label{eq1}
p_{\mu}(y_{[2]}|y_{[1]})=\frac{p_{\mu}(y_{[1]},
y_{[2]})}{p_{\mu}(y_{[1]})}
\end{equation}
the conditional p.m.f. of $y_{[2]}$ given $y_{[1]}$, and $\mathrm{
P}_{\mu}$
the probability distribution, respectively, when $\mu$ is the true
parameter. For any $\varepsilon>0$, we have
%
%e2 #&#
\begin{eqnarray}\label{eq2}
\mathrm{ P}_{\mu}\bigl\{p_{\mu}(y_{[1]},y_{[2]})
\leq p_{\mu+\varepsilon}(y_{[1]}, y_{[2]})|y_{[1]}\bigr
\}&=&\mathrm{ P}_{\mu} \biggl\{\frac{p_{\mu
+\varepsilon}(y_{[1]},y_{[2]})}{p_{\mu}(y_{[1]},y_{[2]})}\geq1\Big\vert
y_{[1]} \biggr\}
\nonumber
\\
&\leq&\mathrm{ E} \biggl\{\frac{p_{\mu+\varepsilon}(y_{[1]},
y_{[2]})}{p_{\mu}(y_{[1]},y_{[2]})}\Big\vert y_{[1]}
\biggr\}
\nonumber\\
&=&\sum_{y_{[2]}}\frac{p_{\mu+\varepsilon}(y_{[1]},y_{[2]})}{p_{\mu}(y_{[1]},
y_{[2]})}p_{\mu}(y_{[2]}|y_{[1]})
\\
&=&\sum_{y_{[2]}}\frac{p_{\mu+\varepsilon}(y_{[1]},y_{[2]})}{p_{\mu}(y_{[1]})}
\nonumber
\\
&=&\frac{p_{\mu+\varepsilon}(y_{[1]})}{p_{\mu}(y_{[1]})},\nonumber
\end{eqnarray}
using (\ref{eq1}). A more general form of (\ref{eq2}) is given in Section~\ref{sec5.1}.

On the other hand, by the standard asymptotic arguments [e.g., \citet{Jia10}, page~9],
it can be shown that the likelihood ratio $p_{\mu
+\varepsilon}(y_{[1]})/p_{\mu}(y_{[1]})$ converges to zero in probability,
as $m\wedge n\rightarrow\infty$. Here we use the fact that the components
of $y_{[1]}$, $y_{ii}, 1\leq i\leq m\wedge n$ are independent Bernoulli
random variables. It follows that, for any $\eta>0$, there is $N_{\eta}
\geq1$ such that, with probability $\geq1-\eta$, we have $\zeta_{N}=
\mathrm{ P}_{\mu}\{p_{\mu}(y_{[1]},y_{[2]})\leq p_{\mu+\varepsilon}(y_{[1]},
y_{[2]})|y_{[1]}\}\leq\gamma^{m\wedge n}$ for some $0<\gamma<1$, if
$m\wedge n\geq N_{\eta}$. The argument shows that $\zeta_{N}=O_\mathrm{
P}(\gamma^{m\wedge n})$, hence converges to $0$ in probability. It
follows, by the dominated convergence theorem, that $\mathrm{ E}_{\mu}(
\zeta_{N})=\mathrm{ P}_{\mu}\{p_{\mu}(y_{[1]},y_{[2]})\leq p_{\mu
+\varepsilon}(y_{[1]},y_{[2]})\}\rightarrow0$. Similarly, we have
$\mathrm{
P}_{\mu}\{p_{\mu}(y_{[1]},y_{[2]})\leq p_{\mu-\varepsilon}(y_{[1]},
y_{[2]})\}\rightarrow0$. The rest of the proof follows by the standard
arguments [e.g., \citet{Jia10}, pages~9--10].~%
\end{pf}

The result of Theorem~\ref{th1} is usually referred to as Cram\'{e}r-type
consistency [Cram\'{e}r (\citeyear{Cra46})], which states that a root to the likelihood
equation is consistent. However, it does not always imply that the MLE,
which by definition is the (global) maximizer of the likelihood function,
is consistent. A stronger result is called Wald-type consistency
[\citet{Wal49}; also see \citet{Wol49}], which states that the MLE is consistent.
Note that the limiting process in Theorem~\ref{th1} is $m, n\rightarrow\infty$,
or, equivalently, $m\wedge n\rightarrow\infty$ (see Section~\ref{sec5.4} for
discussion). With a slightly more restrictive limiting process, the
Wald-consistency can actually be established, as follows.

\begin{theorem}[(Wald consistency)]\label{th2}
If $(m\wedge n)^{-1}\log(m\vee n)
\rightarrow0$ as $m,n\rightarrow\infty$, then the MLE of $\mu$ is
consistent.
\end{theorem}

\begin{pf}
Define $p_{0}(\lambda)=\mathrm{ E}\{h(\lambda+\xi)\}$, where
$h(x)=e^{x}/(1+e^{x})$ and $\xi\sim N(0,2)$. Write $p_{0}=p_{0}(\mu)$.
For any integer $k$, divide the interval $[k,k+1)$ by
$\lambda_{k,j}=k+\delta(mn)^{-1}(m\wedge n)j$, $j=1,\ldots,J$, where $J=
[mn/\delta(m\wedge n)]$ and $0<\delta<1-p_{0}$. It is easy to show that
$|(\partial/\partial\mu)\log p_{\mu}(y_{[1]},y_{[2]})|\leq mn$ uniformly
for all $\mu$. Thus, for any $\lambda\in[k,k+1)$, there is $1\leq j\leq
J$, such that $\log p_{\lambda}(y_{[1]},y_{[2]})-\log p_{\lambda_{k,j}}(
y_{[1]},y_{[2]})=\{(\partial/\partial\mu)\log p_{\mu}(y_{[1]},y_{[2]})
|_{\mu=\tilde{\lambda}}\}(\lambda-\lambda_{k,j})\leq\delta(m\wedge n)$,
where $\tilde{\lambda}$ lies between $\lambda$ and $\lambda_{k,j}$. It
follows that
\[
\sup_{\lambda\in[k,k+1)}\frac{p_{\lambda}(y_{[1]},y_{[2]})}{p_{\mu}(y_{[1]},
y_{[2]})}\leq e^{\delta(m\wedge n)}\max_{1\leq j\leq J}
\frac{p_{\lambda_{k,
j}}(y_{[1]},y_{[2]})}{p_{\mu}(y_{[1]},y_{[2]})}.
\]
Therefore, by the subset argument [see (\ref{eq2})], we have
%
%e3 #&#
\begin{eqnarray}\label{eq3}
&&\mathrm{ P}_{\mu} \biggl\{\sup_{\lambda\in[k,k+1)}
\frac{p_{\lambda}(
y_{[1]},y_{[2]})}{p_{\mu}(y_{[1]},y_{[2]})}>1\Big\vert y_{[1]} \biggr\}
\nonumber
\\
&&\qquad\leq\sum_{j=1}^{J}\mathrm{
P}_{\mu} \biggl\{\frac{p_{\lambda_{k,j}}(
y_{[1]},y_{[2]})}{p_{\mu}(y_{[1]},y_{[2]})}>e^{-\delta(m\wedge n)}\Big\vert
y_{[1]} \biggr\}
\\
&&\qquad \leq e^{\delta(m\wedge n)}\sum_{j=1}^{J}
\frac{p_{\lambda_{k,
j}}(y_{[1]})}{p_{\mu}(y_{[1]})}.\nonumber
\end{eqnarray}

On the other hand, we have $0\leq1-p_{0}(\lambda)=\mathrm{ E}\{1+\exp
(\lambda
+\xi)\}^{-1}\leq \break e^{-\lambda}\mathrm{ E}(e^{-\xi})=e^{1-\lambda}$; and,
similarly, $0\leq p_{0}(\lambda)\leq e^{1+\lambda}$. Let ${\cal
A}_{\delta}
=\{|\Delta|\leq\delta\}$ with $\Delta=(m\wedge n)^{-1}\sum_{i=1}^{m\wedge
n}y_{ii}-p_{0}$. If $k \geq1$, then, for any $1\leq j\leq J$, write $p_{1}
=p_{0}(\lambda_{k,j})$. We have, on ${\cal A}_{\delta}$,
\begin{eqnarray*}
\frac{p_{\lambda_{k,j}}(y_{[1]})}{p_{\mu}(y_{[1]})}&=& \biggl\{ \biggl(\frac
{p_{1}}{p_{0}} \biggr)^{p_{0}+\Delta}
\biggl(\frac{1-p_{1}}{1-p_{0}} \biggr)^{1-p_{0}-\Delta} \biggr\}^{m\wedge n}
\\
&\leq&\bigl\{a_{\delta}^{-1}(1-p_{1})^{1-p_{0}-\delta}
\bigr\}^{m\wedge n}
\\
&\leq&\bigl[a_{\delta}^{-1}\exp\bigl\{(1-\lambda_{k,j})
(1-p_{0}-\delta)\bigr\} \bigr]^{m\wedge n}
\\
&\leq&\exp \bigl[\bigl\{1-p_{0}-\delta-\log a_{\delta}-(1-p_{0}-
\delta)k\bigr\} (m\wedge n) \bigr],
\end{eqnarray*}
where $a_{\delta}=\inf_{|x|\leq\delta}p_{0}^{p_{0}+x}(1-p_{0})^{1
-p_{0}-x}>0$. It follows, by (\ref{eq3}), that
\begin{eqnarray*}
&&\mathrm{ P}_{\mu} \biggl\{\sup_{\lambda\in[k,k+1)}
\frac{p_{\lambda}(
y_{[1]},y_{[2]})}{p_{\mu}(y_{[1]},y_{[2]})}>1\Big\vert y_{[1]} \biggr\}
\\
&&\qquad \leq\frac{mn}{\delta(m\wedge n)}\exp\bigl[\bigl\{1-p_{0}-\log
a_{\delta}-(1-p_{0} -\delta)k\bigr\}(m\wedge n)\bigr]
\end{eqnarray*}
on ${\cal A}_{\delta}$, or, equivalently, that
%
%e4 #&#
\begin{eqnarray}\label{eq4}
&&\mathrm{ P}_{\mu} \biggl\{\sup_{\lambda\in[k,k+1)}
\frac{p_{\lambda}(
y_{[1]},y_{[2]})}{p_{\mu}(y_{[1]},y_{[2]})}>1, |\Delta|\leq\delta\Big\vert y_{[1]} \biggr\}
\nonumber
\\[-8pt]
\\[-8pt]
\nonumber
&&\qquad \leq\frac{mn}{\delta(m\wedge n)}\exp\bigl[\bigl\{1-p_{0}-\log
a_{\delta}-(1-p_{0} -\delta)k\bigr\}(m\wedge n)
\bigr]1_{{\cal A}_{\delta}}.
\end{eqnarray}
Note that ${\cal A}_{\delta}\in{\cal F}(y_{[1]})$. By taking expectations
on both sides of (\ref{eq4}), it follows that the unconditional probability
corresponding to the left side is bounded by the right side without
$1_{{\cal A}_{\delta}}$, for $k=1,2,\ldots.$ Therefore, we have
%
%e5 #&#
\begin{eqnarray}\label{eq5}
&&\mathrm{ P}_{\mu} \biggl\{\sup_{\lambda\in[k,k+1)}\frac{p_{\lambda}(y_{[1]},
y_{[2]})}{p_{\mu}(y_{[1]},y_{[2]})}>1
 \mbox{ for some } k\geq K, |\Delta|\leq\delta \biggr\}
\nonumber
\\
&&\qquad\leq\sum_{k=K}^{\infty}\mathrm{
P}_{\mu} \biggl\{\sup_{\lambda\in
[k,k+1)}\frac
{p_{\lambda}(y_{[1]},y_{[2]})}{p_{\mu}(y_{[1]},y_{[2]})}>1,|\Delta|\leq
\delta \biggr\}
\nonumber
\\
&&\qquad\leq\frac{mn}{\delta(m\wedge n)}\exp\bigl\{(1-p_{0}-\log a_{\delta
}) (m
\wedge n)\bigr\}\sum_{k=K}^{\infty}e^{-(1-p_{0}-\delta)(m\wedge n)k}
\nonumber\\
&&\qquad=\frac{mn}{\delta(m\wedge n)}\exp\bigl\{(1-p_{0}-\log a_{\delta}) (m
\wedge n)\bigr\}\frac{e^{-(1-p_{0}-\delta)(m\wedge n)K}}{1-e^{-(1-p_{0}-\delta
)(m\wedge
n)}}\\
&&\qquad= \bigl\{1-e^{-(1-p_{0}-\delta)(m\wedge n)} \bigr\}^{-1}\nonumber\\
&&\qquad\quad{}\times\exp\bigl[ -(m\wedge n)\bigl
\{(1-p_{0}-\delta)K-1+p_{0}+\log a_{\delta}\nonumber\\
&&\hspace*{80pt}\qquad\quad{}-(m\wedge n)^{-1}\log(m\vee
n)+(m\wedge n)^{-1}\log\delta\bigr\}\bigr].\nonumber
\end{eqnarray}
Thus, if we choose $K$ such that $(1-p_{0}-\delta)K-1+p_{0}+\log
a_{\delta}\geq1$, then, for large $m\wedge n$, the probability on the
left side of (\ref{eq5}) is bounded by $2e^{-(m\wedge n)/2}$. On the other
hand, we have $\mathrm{ P}_{\mu}({\cal A}_{\delta}^{c})\rightarrow0$, as
$m\wedge n\rightarrow\infty$. Thus, we have
%
%e6 #&#
\begin{eqnarray}\label{eq6}
&&\mathrm{ P} \biggl\{\frac{p_{\lambda}(y_{[1]},y_{[2]})}{p_{\mu}(y_{[1]},
y_{[2]})}>1 \mbox{ for some } \lambda
\geq K \biggr\}
\nonumber
\\[-8pt]
\\[-8pt]
\nonumber
&&\qquad\leq 2e^{-(m\wedge n)/2}+\mathrm{ P}\bigl({\cal A}_{\delta
}^{c}
\bigr)\longrightarrow0
\end{eqnarray}
as $m\wedge n\rightarrow\infty$. Similarly, the left side of (\ref{eq6}), with the
words ``$\lambda\geq K$'' replaced by ``$\lambda\leq-K$,'' goes to zero,
as $m\wedge n\rightarrow\infty$, if $K$ is chosen sufficiently large.

On the other hand, again by the subset argument, it can be shown (see
the supplementary material [\citet{supp}]) that for any $\varepsilon>0$ and $K>|\mu
|+\varepsilon$,
we have
%
%e7 #&#
\begin{eqnarray}\label{eq7}
P_{\mu} \biggl\{\sup_{\lambda\in[-K,\mu-\varepsilon)\cup(\mu+\varepsilon,K]} \frac{p_{\lambda}(y_{[1]},y_{[2]})}{p_{\mu}(y_{[1]},y_{[2]})}>1 \biggr\} &
\longrightarrow&0
\end{eqnarray}
as $m, n\rightarrow\infty$. The consistency of the MLE then follows by
combining (\ref{eq7}) with the previously proved results.
\end{pf}

%s3 #&#
\section{Beyond}\label{sec3}
We consider a few more applications of the subset argument, introduced
in the previous section. All applications are regarding a general GLMM,
whose definition is given below for the sake of completeness [see, e.g., \citet{Jia07} for further details].

(i) Suppose that, given a vector $u$ of random effects, responses
$y_{1},\ldots,y_{N}$ are conditionally independent with conditional
density function,\vspace*{1pt} with respect to a $\sigma$-finite measure $\nu$,
given by the exponential family $f_{i}(y_{i}|u)=\exp[a_{i}^{-1}(
\phi)\{y_{i}\xi_{i}-b(\xi_{i})\}+c_{i}(y_{i},\phi)]$, where $\phi$
is a dispersion parameter (which in some cases is known), and $b(\cdot),
a_{i}(\cdot), c_{i}(\cdot,\cdot)$ are known, continuously differentiable
functions with respect to $\xi_{i}$ and $\phi$. The natural parameter
of the conditional exponential family, $\xi_{i}$, is therefore
associated with the conditional mean, $\mu_{i}=\mathrm{ E}(y_{i}|u)$,
according to the properties of the exponential family [e.g., \citet{McCNel83}, Section 2.2.2]. (ii) Furthermore, suppose that $\mu_{i}$
satisfies $g(\mu_{i})=x_{i}'\beta+z_{i}'u$, where $x_{i}, z_{i}$ are
known vectors, $\beta$ is a vector of unknown parameters, and $g(\cdot)$
is a link function. (iii) Finally, assume that $u\sim N(0,G)$, where the
covariance matrix $G$ may depend on a vector $\varphi$ of dispersion
parameters.

It is typically possible to find a subset of the data that are independent,
in some way, under a general GLMM. For example, under the so-called ANOVA
GLMM [e.g., \citet{Lin97}], a subset of independent data can always be found.
Here an ANOVA GLMM satisfies $g(\mu)=X\beta+Z_{1}u_{1}+\cdots+Z_{s}u_{s}$,
where $\mu=(\mu_{i})_{1\leq i\leq N}$, $g(\mu)=[g(\mu_{i})]_{1\leq i\leq
N}$, $X=(x_{i}')_{1\leq i\leq N}$, $Z_{r}=(z_{ir}')_{1\leq i\leq N},
1\leq
r\leq s$, are known matrices, $u_{r}, 1\leq r\leq s$ are vectors of
independent random effects, and $u_{1},\ldots,u_{s}$ are independent.
Examples~\ref{ex1} and~\ref{ex2} are special cases of the ANOVA GLMM. Note that in both
examples the responses are indexed by $(i,j)$, instead of $i$, but this
difference is trivial. Nevertheless, the ``trick'' is to select a subset,
or more than one subsets if necessary, with the following desirable
properties: (I) the subset(s) can be divided into independent
clusters with the number(s) of clusters increasing with the
sample size; and (II) the combination of the subset(s) jointly
identify all the unknown parameters. More specifically, let $y_{i}^{(a)},
i=1,\ldots,N_{a}$ be the $a$th subset of the data, $1\leq a \leq b$,
where $b$ is a fixed positive integer. Suppose that, for each $a$,
there is a partition, $\{1,\ldots,N_{a}\}=\bigcup_{j=1}^{m_{a}}S_{a,j}$.
Let $y_{a,j}=[y_{i}^{(a)}]_{i\in S_{a,j}}$, and $p_{\theta}(y_{a,j})$ be
the probability density function (p.d.f.) of $y_{a,j}$, with respect to the
measure $\nu$ (or the product measure induced by $\nu$ if $y_{a,j}$ is
multivariate), when $\theta$ is the true parameter vector. Let $\Theta$
denote the parameter space, and $\theta_{0}$ the true parameter vector.
Then, (I) and (II) can be formally stated as follows:

\begin{longlist}[(A1)]
\item[(A1)] $y_{a,j}, 1\leq j\leq m_{a}$ are independent with $m_{a}
\rightarrow\infty$ as $N\rightarrow\infty, 1\leq a\leq b$;

\item[(A2)] for every $\theta\in\Theta\setminus\{\theta_{0}\}$, we have
\[
\min_{1\leq a\leq b}\limsup_{N\rightarrow\infty}\frac{1}{m_{a}} \sum
_{j=1}^{m_{a}}\mathrm{ E}_{\theta_{0}} \biggl[\log
\biggl\{\frac
{p_{\theta}(y_{a,j})}{p_{\theta_{0}}(y_{a,j})} \biggr\} \biggr]< 0.
\]
\end{longlist}

Note that (A2) controls the average Kullback--Leibler information
[\citet{KulLei51}]; thus, the inequality always holds if $<$ is
replaced by $\leq$.
%s3.1 #&#
\subsection{Finite parameter space}\label{sec3.1}

Let us first consider a simpler case by assuming that $\Theta$ is finite.
Although the assumption may seem restrictive, it is not totally
unrealistic. For example, any computer system only allows a finite number
of digits. This means that the parameter space that is practically stored
in a computer system is finite. Using the subset argument, it is fairly
straightforward to prove the following (see the supplementary material [\citet{supp}]).

\begin{theorem}\label{th3}
Under assumptions \textup{(A1)} and \textup{(A2)}, if, in
addition,
\begin{longlist}[(A3)]
\item[(A3)] for every $\theta\in\Theta\setminus\{\theta_{0}\}$, we have
\[
\frac{1}{m_{a}^{2}}\sum_{j=1}^{m_{a}}
\operatorname{ var}_{\theta
_{0}} \biggl[ \log \biggl\{\frac{p_{\theta}(y_{a,j})}{p_{\theta_{0}}(y_{a,j})} \biggr
\} \biggr]\longrightarrow 0, \qquad 1\leq a\leq b,
\]
then $\mathrm{ P}_{\theta_{0}}(\hat{\theta}=\theta_{0})\rightarrow1$, as
$N\rightarrow\infty$, where $\hat{\theta}$ is the MLE of $\theta$.
\end{longlist}
\end{theorem}
%s3.2 #&#
\subsection{Euclidean parameter space}\label{sec3.2}

We now consider the case that $\Theta$ is a convex subspace of $R^{d}$,
the $d$-dimensional Euclidean space, in the sense that $\theta_{1},
\theta_{2}\in\Theta$ implies $(1-t)\theta_{1}+t\theta_{2}\in\Theta$ for
every $t\in(0,1)$. In this case, we need to strengthen assumptions
(A2), (A3) to the following:

\begin{longlist}[(B3)]
\item[(B2)] $\theta_{0}\in\Theta^\mathrm{ o}$, the interior of $\Theta$,
and there is $0<M<\infty$ [same as in (B3) below] such that, for
every $\varepsilon>0$, we have
%
%e8 #&#
\begin{eqnarray}\label{eq8}
\limsup_{N\rightarrow\infty}\sup_{\theta\in\Theta,\varepsilon\leq
|\theta-\theta_{0}|\leq M}\min_{1\leq a\leq b}\frac{1}{m_{a}}
\sum_{j=1}^{m_{a}}\mathrm{ E}_{\theta_{0}}
\biggl[\log \biggl\{\frac
{p_{\theta}(y_{a,j})}{p_{\theta_{0}}(y_{a,j})} \biggr\} \biggr]<0.
\end{eqnarray}

\item[(B3)] There are positive constant sequences $s_{N}, s_{a,N},
1\leq
a\leq b$ such that
%
%e9 #&#
\begin{eqnarray}\label{eq9}
\sup_{\theta\in\Theta,|\theta-\theta_{0}|\leq M}\max_{1\leq c\leq d} \biggl\llvert \frac{\partial}{\partial\theta_{c}}
\log\bigl\{p_{\theta}(y)\bigr\}\biggr\rrvert = O_\mathrm{
P}(s_{N})
\end{eqnarray}
with $\log(s_{N})/\min_{1\leq a\leq b}m_{a}\rightarrow0$,
where $p_{\theta}(y)$ is the p.d.f. of $y=(y_{i})_{1\leq i\leq N}$ given
that $\theta=(\theta_{c})_{1\leq c\leq d}$ is the true parameter
vector,
%
%e10 #&#
\begin{eqnarray}\label{eq10}
\sup_{\theta\in\Theta,|\theta-\theta_{0}|\leq M}\frac
{1}{m_{a}}\sum_{j=1}^{m_{a}}
\max_{1\leq c\leq d}\biggl\llvert \frac
{\partial}{\partial\theta_{c}}\log\bigl
\{p_{\theta}(y_{a,j})\bigr\}\biggr\rrvert = o_\mathrm{
P}(s_{a,N})
\end{eqnarray}
with $\log(s_{a,N})/m_{a}\rightarrow0$; and (for the same $s_{a,N}$)
%
%e11 #&#
\begin{eqnarray}\label{eq11}\qquad
\sup_{\theta\in\Theta,|\theta-\theta_{0}|\leq M}\frac{s_{a,N}^{d
-1}}{m_{a}^{2}}\sum_{j=1}^{m_{a}}
\operatorname{ var}_{\theta_{0}} \biggl[\log \biggl\{\frac{p_{\theta}(y_{a,j})}{p_{\theta_{0}}(y_{a,j})} \biggr\}
\biggr] \longrightarrow0,\qquad  1\leq a\leq b.
\end{eqnarray}
\end{longlist}

\begin{theorem}\label{th4}
Under assumptions \textup{(A1), (B2)} and
\textup{(B3)},
there is, with probability $\rightarrow1$, a root to the likelihood
equation, $\hat{\theta}$, such that $\hat{\theta}\stackrel{\mathrm{
P}}{\longrightarrow}\theta_{0}$, as $N\rightarrow\infty$.
\end{theorem}

\begin{pf}
Aside from the use of the subset argument, the lines of the
proof are similar to, for example, the standard arguments of Lehmann and
Casella [(\citeyear{LehCas98}), the beginning part of the proof of Theorem 5.1], although
some details are more similar to \citet{Wol49}. We outline the key steps
below and refer the details to the supplementary material [\citet{supp}]. Once again, the
innovative part is the consideration of the conditional probability given
the subset data and, most importantly, the subset inequality (\ref{eq15}) in the
sequel.

For any $\varepsilon>0$, assume, without loss of generality,
that $\{\theta\dvtx |\theta-\theta_{0}|\leq\varepsilon\}\subset\Theta$ and
$C_{\varepsilon}=\{\theta\in R^{d}\dvtx|\theta_{c}-\theta_{0c}|\leq
\varepsilon, 1\leq c\leq d\}\subset\{\theta\in\Theta\dvtx |\theta
-\theta_{0}|\leq M\}$. Essentially, all we need to show is that, as
$N\rightarrow\infty$,
%
%e12 #&#
\begin{equation}\label{eq12}
P(\varepsilon)\equiv\mathrm{ P}_{\theta_{0}} \Bigl\{p_{\theta_{0}}(y)\leq
\sup_{\theta\in\partial C_{\varepsilon}}p_{\theta}(y) \Bigr\} \longrightarrow 0,
\end{equation}
where $\partial C_{\varepsilon}$ is the boundary of $C_{\varepsilon}$, which
consists of $\theta\in C_{\varepsilon}$ such that $|\theta_{c}-\theta_{0c}|
=\varepsilon$ for some $1\leq c\leq d$. Define
\[
S_{N,a}(\theta)=\frac{1}{m_{a}}\sum_{j=1}^{m_{a}}
\mathrm{ E}_{\theta_{0}} \biggl[\log \biggl\{\frac{p_{\theta}(y_{a,
j})}{p_{\theta_{0}}(y_{a,j})} \biggr\}
\biggr], \qquad 1\leq a\leq b,
\]
and $I_{N}(\theta)=\min\{1\leq a\leq b\dvtx S_{N,a}(\theta)=\min_{1\leq
a'\leq
b}S_{N,a'}(\theta)\}$. Then, $\partial C_{\varepsilon}=\break\bigcup_{a=1}^{b}\partial
C_{\varepsilon}\cap\Theta_{N,a}$, where $\Theta_{N,a}=\{\theta\in\Theta:
I_{N}(\theta)=a\}$. Then, we have
%
%e13 #&#
\begin{equation}\label{eq13}
P(\varepsilon)\leq\sum_{a=1}^{b}
\mathrm{ P}_{\theta_{0}} \Bigl\{ p_{\theta_{0}}(y) \leq\sup_{\theta\in\partial C_{\varepsilon}\cap\Theta_{N,a}}p_{\theta}(y)
\Bigr\}.
\end{equation}

For a fixed $1\leq a\leq b$, let $\delta$ be a small, positive number
to be determined latter, and $K=[e^{\delta m_{a}}]+1$. For any $l=(l_{1},
\ldots,l_{d})$, where $0\leq l_{c}\leq K-1, 1\leq c\leq d$, select a point
$\theta_{l}$ from the subset $\{\theta\dvtx \theta_{0c}-\varepsilon
+2\varepsilon
l_{c}/K\leq\theta_{c}\leq\theta_{0c}-\varepsilon+2\varepsilon
(l_{c}+1)/K, 1\leq
c\leq d\}\cap\partial C_{\varepsilon}\cap\Theta_{N,a}$, if the latter
is not
empty; otherwise, do not select. Let $D$ denote the collection of all such
points. Also let $B$ denote the left side of~(\ref{eq9}). It can be shown that
%
%e14 #&#
\begin{eqnarray}\label{eq14}
&&\mathrm{ P}_{\theta_{0}} \Bigl\{p_{\theta_{0}}(y)\leq\sup_{\theta\in
\partial
C_{\varepsilon}\cap\Theta_{N,a}}p_{\theta}(y)
\Bigr\}
\nonumber
\\[-8pt]
\\[-8pt]
\nonumber
&&\qquad\leq\mathrm{ P}_{\theta_{0}} \biggl\{\exp \biggl(\frac{2d\varepsilon
B}{K} \biggr)
>2 \biggr\}+\mathrm{ P}_{\theta_{0}} \Bigl\{p_{\theta_{0}}(y)\leq2
\max_{\theta
\in D}p_{\theta}(y) \Bigr\}.
\end{eqnarray}

We now apply the subset argument. Let $y_{[1]}$ denote the combined vector
of $y_{a,j}, 1\leq j\leq m_{a}$, and $y_{[2]}$ the vector of the rest of
$y_{1},\ldots,y_{N}$. Then,\vadjust{\goodbreak} similar to the argument of (\ref{eq2}), we have, for
any $\theta\in D$,
%
%e15 #&#
\begin{equation}\label{eq15}
\mathrm{ P}_{\theta_{0}}\bigl\{p_{\theta_{0}}(y)\leq2p_{\theta}(y)|y_{[1]}
\bigr\} \leq 2\frac{p_{\theta}(y_{[1]})}{p_{\theta_{0}}(y_{[1]})}.
\end{equation}
Using this result, it can be shown that $\mathrm{ P}_{\theta_{0}}\{
p_{\theta_{0}}(y)\leq2\max_{\theta\in D}p_{\theta}(y)|y_{[1]}\}
=o_\mathrm{
P}(1)$. From here, (\ref{eq12}) can be established.
\end{pf}

Again, Theorem~\ref{th4} is a Cram\'{e}r-consistency result. On the other hand,
Wald-consistency can be established under additional assumptions that
control the behavior of the likelihood function in a neighborhood of
infinity. For example, the following result may be viewed as an
extension of Theorem~\ref{th2}. The proof is given in the supplementary material [\citet{supp}].
Once again, the subset argument plays a critical role in the proof.
For simplicity, we focus on the case of discrete responses, which is
typical for GLMMs. In addition, we assume the following. For any
$0\leq v<w$, define $S_{d}[v,w)=\{x\in R^{d}\dvtx v\leq|x|<w\}$ and write,
in short, $S_{d}(k)=S_{d}[k,k+1)$ for $k=1,2,\ldots.$

\begin{longlist}[(C1)]
\item[(C1)] There are sequences of constants, $b_{k}, c_{N}\geq1$, and
random variables, $\zeta_{N}$, where $c_{N}, \zeta_{N}$ do not depend on
$k$, such that $\zeta_{N}=O_\mathrm{ P}(1)$ and
\begin{eqnarray*}
\sup_{\theta\in\Theta\cap S_{d}[k-1,k+2)}\max_{1\leq
c\leq d}\biggl\llvert \frac{\partial}{\partial\theta_{c}}\log
\bigl\{p_{\theta}(y)\bigr\} \biggr\rrvert \leq b_{k}c_{N}
\zeta_{N},\qquad  k=1,2,\ldots
\end{eqnarray*}

\item[(C2)] There is a subset of independent data vectors, $y_{(j)},
1\leq j\leq m_{N}$ [not necessarily among those in (A1)] so that:
(i) $\mathrm{ E}_{\theta_{0}}|\log\{p_{j,\theta_{0}}(y_{(j)})\}|$ is bounded,
$p_{j,\theta}(\cdot)$ being the p.m.f. of $y_{(j)}$ under $\theta$; (ii)
there is a sequence of positive constants, $\gamma_{k}$, with
$\lim_{k\rightarrow\infty}\gamma_{k}=\infty$, and a subset ${\cal
T}_{N}$ of possible values of $y_{(j)}$, such that for every $k\geq1$
and $\theta\in\Theta\cap S_{d}(k)$, there is $t\in{\cal T}_{N}$ satisfying
$\max_{1\leq j\leq m_{N}}\log\{p_{j,\theta}(t)\}\leq-\gamma_{k}$;
(iii) $\inf_{t\in{\cal T}_{N}}m_{N}^{-1}\sum_{j=1}^{m_{N}}p_{j,\theta_{0}}
(t)\geq\rho$ for some constant $\rho>0$; and (iv) $|{\cal T}_{N}|/m_{N}=
o(1)$, and $c_{N}^{d}\sum_{k=K}^{\infty}k^{d_{1}}b_{k}^{d}e^{-\delta
m_{N}\gamma_{k}}=o(1)$ for some $K\geq1$ and $\delta<\rho$, where $d_{1}=
d1_{(d>1)}$.
\end{longlist}

It is easy to verify that the new assumptions (C1), (C2) are
satisfied in the case of Theorem~\ref{th2} for the open problem (see
the supplementary material [\citet{supp}]). Another example is considered in the next section.

\begin{theorem}\label{th5}
Suppose that \textup{(A1)} holds; \textup{(B2), (B3)} hold for
any fixed $M>0$ (instead of some $M>0$), and with the $s_{a,N}^{d-1}$ in
(\ref{eq11}) replaced by $s_{a,N}^{d}$. In addition, suppose that \textup{(C1)},
\textup{(C2)} hold. Then, the MLE of $\theta_{0}$ is consistent.
\end{theorem}

%s4 #&#
\section{Example}\label{sec4}

Let us consider a special case of Example~\ref{ex1} with $x_{ijk}'\beta=\mu$,
but $\sigma^{2}$ and $\tau^{2}$ unknown. We change the notation slightly,
namely, $y_{i,j,k}$ instead of $y_{ijk}$. Suppose that\vadjust{\goodbreak} $S=S_{1}\cup S_{2}$
such that $c_{ij}=r, (i,j)\in S_{r}$, $r=1,2$ (as in the case of the
salamander data). We use two subsets to jointly identify all the
unknown parameters. The first subset is similar to that used in the
proofs of Theorems~\ref{th1} and~\ref{th2}, namely, $y_{i,i}=(y_{i,i,k})_{k=1,2}, (i,
i)\in S_{2}$. Let $m_{1}$ be the total number of such $(i,i)$'s, and
assume that $m_{1}\rightarrow\infty$, as $m,n\rightarrow\infty$. Then,
the subset satisfies (A1). Let $\theta=(\mu,\sigma^{2},\tau^{2})'$.
It can be shown that the sequence $y_{i,i}, (i,i)\in S_{2}$ is a sequence
of i.i.d. random vectors with the probability distribution, under
$\theta$,
given by
%
%e16 #&#
\begin{equation}\label{eq16}
p_{\theta}(y_{i,i})=\mathrm{ E} \biggl[\frac{\exp\{y_{i,i,\cdot}(\mu
+\xi)\}}{\{1+\exp(\mu+\xi)\}^{2}}
\biggr],
\end{equation}
where $\xi\sim N(0,\psi^{2})$, with $\psi^{2}=\sigma^{2}+\tau^{2}$, and
$y_{i,i,\cdot}=y_{i,i,1}+y_{i,i,2}$. By the strict concavity of the
logarithm, we have
%
%e17 #&#
\begin{equation}\label{eq17}
\mathrm{ E}_{\theta_{0}} \biggl[\log \biggl\{\frac{p_{\theta}(y_{i,
i})}{p_{\theta_{0}}(y_{i,i})} \biggr\}
\biggr]< 0
\end{equation}
unless $p_{\theta}(y_{i,i})/p_{\theta_{0}}(y_{i,i})$ is a.s. $\mathrm{
P}_{\theta_{0}}$ a constant, which must be one because both $p_{\theta}$
and $p_{\theta_{0}}$ are probability distributions. It is easy to
show that the probability distribution of (\ref{eq16}) is completely determined
by the function $M(\vartheta)=[M_{r}(\vartheta)]_{r=1,2}$, where
$M_{r}(\vartheta)=\mathrm{ E}\{h_{\vartheta}^{r}(\zeta)\}$ with
$\vartheta=(\mu,
\psi)'$, $h_{\vartheta}(\zeta)=\exp(\mu+\psi\zeta)/\{1+\exp(\mu+\psi
\zeta)
\}$, and $\zeta\sim N(0,1)$. In other words, $p_{\theta}(y_{i,i})=
p_{\theta_{0}}(y_{i,i})$ for all values of $y_{i,i}$ if and only if
$M(\vartheta)=M(\vartheta_{0})$. \citet{Jia98} showed that the function
$M(\cdot)$ is injective [also see \citet{Jia07}, page~221]. Thus, (\ref{eq17}) holds
unless $\mu=\mu_{0}$ and $\psi^{2}=\psi_{0}^{2}$.

It remains to deal with a $\theta$ that satisfies $\mu=\mu_{0}$, $\psi^{2}
=\psi_{0}^{2}$, but $\theta\neq\theta_{0}$. For such a $\theta$, we use
the second subset, defined as $y_{i}=(y_{i,2i-1,1},y_{i,2i,1})'$ such
that $(i,2i-1)\in S$ and $(i,2i)\in S$. Let $m_{2}$ be the total number
of all such $i$'s, and assume that $m_{2}\rightarrow\infty$ as $m, n
\rightarrow\infty$. It is easy to see that (A1) is, again, satisfied
for the new subset. Note that any $\theta$ satisfying $\mu=\mu_{0}$ and
$\psi^{2}=\psi_{0}^{2}$ is completely determined by the parameter
$\gamma
=\sigma^{2}/\psi^{2}$. Furthermore, the new subset is a sequence of i.i.d.
random vectors with the probability distribution, under such a $\theta$,
given by
%
%e18 #&#
\begin{equation}\label{eq18}
p_{\gamma}(y_{i})=\mathrm{ E} \biggl[\frac{\exp\{y_{i,2i-1,1}(\mu_{0}+X)\}}{1+
\exp(\mu_{0}+X)}\cdot
\frac{\exp\{y_{i,2i,1}(\mu_{0}+Y)\}}{1+\exp(\mu_{0}+
Y)} \biggr],
\end{equation}
where $(X,Y)$ has the bivariate normal distribution with $\operatorname
{ var}(X)=
\operatorname{ var}(Y)=\psi_{0}^{2}$ and $\operatorname{
cor}(X,Y)=\gamma$. Similar to (\ref{eq17}),
we have
%
%e19 #&#
\begin{equation}\label{eq19}
\mathrm{ E}_{\gamma_{0}} \biggl[\log \biggl\{\frac
{p_{\gamma}(y_{i})}{p_{\gamma_{0}}(y_{i})} \biggr\}
\biggr]< 0
\end{equation}
unless $p_{\gamma}(y_{i})=p_{\gamma_{0}}(y_{i})$ for all values of $y_{i}$.
Consider (\ref{eq18}) with $y_{i}=(1,1)$ and let $\mathrm{ P}_{\gamma}$ denote the
probability distribution of $(X,Y)$ with the correlation coefficient
$\gamma$. By Fubini's theorem, it can be shown that
%
%e20 #&#
\begin{eqnarray}\label{eq20}\qquad
p_{\gamma}(1,1)=\int_{0}^{\infty}\int
_{0}^{\infty}P_{\gamma}\bigl\{X\geq
\operatorname{ logit}(s)-\mu_{0},Y\geq\operatorname{ logit}(t)-
\mu_{0}\bigr\}\,ds\,dt.
\end{eqnarray}
Hereafter, we refer the detailed derivations to the supplementary material [\citet{supp}].
By Slepian's inequality [e.g., \citet{Jia10}, pages~157--158], the integrand on
the right side of (\ref{eq20}) is strictly increasing with $\gamma$, hence so is
the integral. Thus, if $\gamma\neq\gamma_{0}$, at least we have
$p_{\gamma}(1,1)\neq p_{\gamma_{0}}(1,1)$, hence (\ref{eq19}) holds.

In summary, for any $\theta\in\Theta, \theta\neq\theta_{0}$, we must
have either (\ref{eq17}) or (\ref{eq19}) hold. Therefore, by continuity, assumption
(B2) holds, provided that true variances, $\sigma_{0}^{2}, \tau_{0}^{2}$
are positive. Note that, in the current case, the expectations involved
in (B2) do not depend on either $j$ or $N$, the total sample size.

To verify (B3), it can be shown that $|(\partial/\partial\mu)\log
\{
p_{\theta}(y)\}|\leq N$. Furthermore, we have $|(\partial/\partial
\sigma^{2})\log\{p_{\theta}(y)\}|\vee|(\partial/\partial\tau^{2})\log\{
p_{\theta}(y)\}|\leq(A+C+1)N$ in a neighborhood of $\theta_{0}$, ${\cal
N}(\theta_{0})$. Therefore, (\ref{eq9}) holds with $s_{N}=N$.

As for (\ref{eq10}), it is easy to show that the partial derivatives involved
are uniformly bounded for $\theta\in{\cal N}(\theta_{0})$. Thus, (\ref{eq10})
holds for any $s_{a,N}$ such that $s_{a,N}\rightarrow\infty$, $a=1,2$.
Furthermore, the left side of (\ref{eq11}) is bounded by $c_{a}s_{a,N}^{2}/m_{a}$
for some constant $c_{a}>0$, $a=1,2$ (note that $d=3$ in this case). Thus,
for example, we may choose $s_{a,N}=\sqrt{m_{a}/\{1+\log(m_{a})\}}$, $a=1,
2$, to ensure that $\log(s_{a,N})/m_{a}\rightarrow0$, $a=1,2$, and (\ref{eq11})
holds.

In conclusion, all the assumptions of Theorem~\ref{th4} hold provided that
$\sigma_{0}^{2}>0$, $\tau_{0}^{2}>0$, and $(m_{1}\wedge m_{2})^{-1}\log(
N)\rightarrow0$.

Similarly, the conditions of Theorem~\ref{th5} can be verified. Essentially,
what is new is to check assumptions (C1) and (C2). See
the supplementary material [\citet{supp}].
%s5 #&#
\section{Discussion}\label{sec5}
%s5.1 #&#
\subsection{Remark on subset argument}\label{sec5.1}

In proving a number of results, we have demonstrated the usefulness
of the subset argument. In principle, the method allows one to
argue consistency of the MLE in any situation of dependent data, not
necessarily under a GLMM, provided that one can identify some suitable
subset(s) of the data whose asymptotic properties are easier to
handle, such as collections of independent random vectors. The
connection between the full data and subset data is made by the
subset inequality, which, in a more general form, is a consequence
of the martingale property of the likelihood-ratio [e.g., \citet{Jia10},
pages~244--246]: suppose that $Y_{1}$ is a subvector of a random vector
$Y$. Let $p_{\theta}(\cdot)$ and $p_{1,\theta}(\cdot)$ denote the p.d.f.'s
of $Y$ and $Y_{1}$, respectively, with respect to a $\sigma$-finite
measure $\nu$, under the parameter vector $\theta$. For simplicity,
suppose that $p_{\theta_{0}}, p_{1,\theta_{0}}$ are positive a.e. $\nu$,
and $\lambda(\cdot)$ is a positive, measurable function. Then, for any
$\theta$, we have
\begin{eqnarray*}
\mathrm{ P}_{\theta_{0}}\bigl\{p_{\theta_{0}}(Y)\leq\lambda(Y_{1})
p_{\theta}(Y)|Y_{1}\bigr\}\leq\lambda(Y_{1})
\frac{p_{1,\theta}(Y_{1})}{p_{1,
\theta_{0}}(Y_{1})}\qquad\mbox{a.e. } \nu,
\end{eqnarray*}
where $\mathrm{ P}_{\theta_{0}}$ denotes the probability distribution
corresponding to $p_{\theta_{0}}$.

%s5.2 #&#
\subsection{Quantifying the information loss}\label{sec5.2}
On the other hand, the subset argument merely provides a method of proof
for the consistency of the full-data MLE---it by no means suggests the
subset-data MLE as a replacement for the full-data MLE. In fact, there is
an information loss if such a replacement takes place. To quantify the
information loss, assume the regularity conditions for exchanging the
order of differentiation and integration. Then, the Fisher information
matrix based on the full data can be expressed as
\begin{eqnarray*}
I_\mathrm{ f}(\theta)&=&-\mathrm{ E}_{\theta} \biggl\{
\frac{\partial^{2}}{\partial
\theta\,\partial\theta'}\log p_{\theta}(y) \biggr\}
\\
&=&\mathrm{ E}_{\theta} \biggl[ \biggl\{\frac{\partial}{\partial\theta}\log
p_{\theta}(y) \biggr\} \biggl\{\frac{\partial}{\partial\theta}\log p_{\theta}(y)
\biggr\}' \biggr]-\mathrm{ E}_{\theta} \biggl\{
\frac
{1}{p_{\theta}(y)} \frac{\partial^{2}}{\partial\theta\,\partial\theta'}p_{\theta}(y) \biggr\}
\\
&=&I_{\mathrm{ f},1}(\theta)-I_{\mathrm{ f},2}(\theta).
\end{eqnarray*}
Similarly, the information matrix based on the subset data can be expressed
as $I_\mathrm{ s}(\theta)=I_{\mathrm{ s},1}(\theta)-I_{\mathrm{
s},2}(\theta)$, where
$I_{\mathrm{ s},j}(\theta)$ is $I_{\mathrm{ f},j}(\theta)$ with $y$
replaced by
$y_{[1]}$, $j=1,2$ [$p_{\theta}(y_{[1]})$ denotes the p.d.f. (or
p.m.f.) of
$y_{[1]}$]. By conditioning on $y_{[1]}$, it can be shown that
$I_{\mathrm{
f},2}(\theta)=I_{\mathrm{ s},2}(\theta)$, while $I_{\mathrm{
f},1}(\theta)\geq
I_{\mathrm{ s},1}(\theta)$. It follows that
%
%e21 #&#
\begin{equation}\label{eq21}
I_\mathrm{ f}(\theta)\geq I_\mathrm{ s}(\theta)
\end{equation}
for all $\theta$. Here the inequality means that the difference between
the left side and right side is a nonnegative definite matrix. (\ref{eq21}) suggests
that the information contained in the full data is no less than that
contained in the subset data, which, of course, is what one would expect.
Furthermore, the information loss is given by
%
%e22 #&#
\begin{equation}\label{eq22}
I_\mathrm{ f}(\theta)-I_\mathrm{ s}(\theta)=\mathrm{
E}_{\theta} \biggl[\operatorname{ Var}_{\theta} \biggl\{
\frac{\partial}{\partial\theta}\log p_{\theta}(y) \Big\vert y_{[1]} \biggr\}
\biggr],
\end{equation}
where $\operatorname{ Var}_{\theta}(\cdot|y_{[1]})$ denotes the conditional
covariance matrix given $y_{[1]}$ under $\theta$. The derivations of (\ref{eq21})
and (\ref{eq22}) are deferred to the supplementary material [\citet{supp}]. It is seen from (\ref{eq22})
that the information loss is determined by how much (additional) variation
there is in the score function, $(\partial/\partial\theta)\log p_{\theta}(
y)$, given the subset data $y_{[1]}$. In particular, if $y_{[1]}=y$, then
the score function is a constant vector given $y_{[1]}$ (and $\theta$);
hence $\operatorname{ Var}_{\theta}\{(\partial/\partial\theta)\log
p_{\theta}(y)|
y_{[1]}\}=0$, thus, there is no information loss. In general, of course,
the subset data $y_{[1]}$ is not chosen as $y$; therefore, there will be
some loss of information.

Nevertheless, the information contained in the subset data is usually
sufficient for identifying at least some of the parameters. Note that
consistency is a relatively weak asymptotic property in the sense that
various estimators, including those based on the subset data and, for
example, the method of moments estimator of \citet{Jia98}, are consistent,
even though they may not be asymptotically efficient. Essentially, for
the consistency property to hold, one needs that, in spite of the
potential information loss, the remaining information that the
estimator is able to utilize grows with the sample size. For example,
in the open problem (Sections~\ref{sec1} and~\ref{sec2}), the information contained in
$y_{ii}$ grows at the rate of $m\wedge n$, which is sufficient for
identifying $\mu$; in the example of Section~\ref{sec4}, the information contained
in $y_{i,i}$ grows in the order of $m_{1}$, which is sufficient for
identifying $\mu$ and $\psi^{2}=\sigma^{2}+\tau^{2}$, while the information
contained in $y_{i}$ grows at the rate of $m_{2}$, which is sufficient for
identifying $\gamma=\sigma^{2}/\psi^{2}$. The identification of the
``right'' subset in a given problem is usually suggested by the nature of
the parametrization. As mentioned (see the third paragraph of Section~\ref{sec3}),
a subset $y_{[1]}$ of independent data can always be found under the ANOVA
GLMM (e.g., starting with the first observation, $y_{1}$, one finds
the next observation such that it involves different random effects from
those related to $y_{1}$, and so on). If the $y_{[1]}$ is such that
$\liminf_{N\rightarrow\infty}\lambda_{\min}\{I_\mathrm{ s}(\theta)\}
=\infty$,
where $I_\mathrm{ s}(\theta)$ is as in (\ref{eq21}) and $\lambda_{\min}$
denotes the
smallest eigenvalue, the subset $y_{[1]}$ is sufficient for identifying all
the components of $\theta$; otherwise, more than one subsets are needed in
order to identify all the parameters, as is shown in Section~\ref{sec4}.
%s5.3 #&#
\subsection{Note on computation of MLE}\label{sec5.3}

The subset argument offers a powerful tool for establishing consistency
of the MLE in GLMM with crossed random effects. Note that the idea has
not followed the traditional path of attempting to develop a
(computational) procedure to approximate the MLE. In fact, this might
explain why the computational advances over the past two decades [see,
e.g., \citet{Jia07}, Section 4.1 for an overview] had not led
to a major theoretical breakthrough for the MLE in GLMM in terms of
asymptotic properties. Note that the MLE based on the subset data is a
consistent estimator of the true parameter, and in that sense it is an
approximation to the MLE based on the full data (two consistent estimators
of the same parameter approximate each other). However, there is an
information loss, as discussed in the previous subsection [see (\ref{eq22})], so
one definitely wants to do better computationally.

One computational method that has been developed for computing the MLE
in GLMMs, including those with crossed random effects, is Monte Carlo
EM algorithm [e.g., McCullogh (\citeyear{McC94,McC97}), \citet{BooHob99}]. Here,
however, we would like to discuss another, more recent, computational
advance, known as \textit{data cloning} [DC; \citet{LelDenLut07}, \citet{LelNadSch10}].
The DC uses the Bayesian computational approach for frequentist purposes.
Let $\pi$ denote the prior density function of $\theta$. Then, one has
the posterior,
%
%e23 #&#
\begin{equation}\label{eq23}
\pi(\theta|y)=\frac{p_{\theta}(y)\pi(\theta)}{p(y)},
\end{equation}
where $p(y)$ is the integral of the numerator with respect to $\theta$,
which does not depend on $\theta$. There are computational tools
using the Markov chain Monte Carlo for posterior simulation that
generate random variables from the posterior without having to compute
the numerator or denominator of (\ref{eq23}) [e.g., Gilks, Richardson and
  Spiegelhalter (\citeyear{GiRiSp96});
\citet{Spietal04}]. Thus, we can assume that one can
generate random variables from the posterior. If the observations
$y$ were repeated independently from $K$ different individuals such
that all of these individuals result in exactly the same data, $y$,
denoted by $y^{(K)}=(y,\ldots,y)$, then the posterior based on $y^{(K)}$
is given by
%
%e24 #&#
\begin{equation}\label{eq24}
\pi_{K}\bigl\{\theta|y^{(K)}\bigr\}=\frac{\{p_{\theta}(y)\}^{K}\pi(\theta)}{\int
\{
p_{\theta}(y)\}^{K}\pi(\theta)\,d\theta}.
\end{equation}
\citet{LelDenLut07}, \citet{LelNadSch10} showed that, as $K$ increases, the
right side of (\ref{eq24}) converges to a multivariate normal distribution whose
mean vector is equal to the MLE, $\hat{\theta}$, and whose covariance
matrix is approximately equal to $K^{-1}I_\mathrm{ f}^{-1}(\hat{\theta})$.
Therefore, for large $K$, one can approximate the MLE by the sample mean
vector of, say, $\theta^{(1)},\ldots,\theta^{(B)}$ generated from the
posterior distribution (\ref{eq24}). Denoted the sample mean by $\bar
{\theta}^{(\cdot)}$, and call it the DC MLE. Furthermore, $I_\mathrm{ f}^{-1}
(\hat{\theta})$ [see (\ref{eq21}), (\ref{eq22})] can be approximated by $K$ times the
sample covariance matrix of $\theta^{(1)},\ldots,\theta^{(B)}$. \citet{Tor12}
successfully applied the DC method to obtain the MLE for the
salamander-mating data.

Note that the DC MLE is an approximate, rather than exact, MLE, in the
sense that, as $K\rightarrow\infty$, the difference between $\bar
{\theta}^{(\cdot)}$ and the exact MLE vanishes. Because we have
established consistency of the exact MLE, it follows that
the DC MLE is a consistent estimator as long as the number
$K$ increase with the sample size. More precisely, it is shown
in the supplementary material [\citet{supp}] that, for every $\varepsilon,\delta>0$, there is
$N_{\varepsilon,\delta}$ such that for every $n\geq N_{\varepsilon
,\delta}$ and
$B\geq1$, there is $K(n,B)$ such that $\mathrm{ P}\{|\bar{\theta
}^{(\cdot)}
-\theta_{0}|\geq\varepsilon\}<\delta$, if $K\geq K(n,B)$, where $\theta_{0}$
is the true parameter vector. Note that, as far as consistency is
concerned, one does not need that $B$ goes to infinity. This makes sense
because, as $K\rightarrow\infty$, the posterior (\ref{eq24}) is becoming degenerate
[the asymptotic covariance matrix is $K^{-1}I_\mathrm{ f}^{-1}(\hat
{\theta})$];
thus, one does not need a large $B$ to ``average out'' the variation in
$\bar{\theta}^{(\cdot)}$. Thus, from an asymptotic point of view, the
result of the current paper provides a justification for the DC method.

More importantly, because $B,K$ are up to one's choice, one can make sure
that they are large enough so that there is virtually no information loss,
as was concerned earlier. In this regard, a reasonably large $B$ would
reduce the sampling variation and therefore improve the DC approximation,
and make the computation more efficient. See \citet{LelNadSch10} for
discussion on how to choose $B$ and $K$ from practical points of view.\vadjust{\goodbreak}

As for the prior $\pi$, \citet{LelNadSch10} only suggests that it
be chosen according to computational convenience and be proper (to avoid
improper posterior). Following the subset idea, an obvious choice for
the prior would be the multivariate normal distribution with mean vector
$\hat{\theta}_\mathrm{ s}$, the subset-data MLE, and covariance matrix
$I_\mathrm{ s}^{-1}(\hat{\theta}_\mathrm{ s})$ [defined above (\ref{eq21})].
Note that
$I_\mathrm{ s}(\theta)$ is much easier to evaluate than $I_\mathrm{
f}(\theta)$.
This would make the procedure more similar to the empirical Bayes than
the hierarchical one. Nevertheless, the DC only uses the Bayesian
computational tool, as mentioned.

%s5.4 #&#
\subsection{Regarding the limiting process}\label{sec5.4}

In some applications of GLMM, the estimation of the random effects
are of interest. There have also been developments in semiparametric
GLM and nonparametric ANOVA. In those cases, the random effects are
treated the same way as the fixed effects. As a result, the proof of the
consistency results in those cases usually impose constraints on the
ratio of the number of effects and number of observations falling in each
cluster [e.g., \citet{Che95,Jia99,WuLia04}, and Wang, Tsai and Qu (\citeyear{WanTsaQu12})].
A major difference exists, however, between the case of
clustered data (e.g., Example~\ref{ex2}) and that with crossed random effects
(e.g., Example~\ref{ex1}) in that, in the latter case, the data cannot be
divided into independent groups (with the number of groups increasing
with the sample size). Furthermore, the necessary constraints are very
different depending on the interest of estimation. Consider, for example,
a very simple case of linear mixed model, $y_{ij}=\mu+u_{i}+v_{j}+e_{ij}$,
$i=1,\ldots,m, j=1,\ldots,n$, where the $u_{i}$'s and $v_{j}$'s are random
effects, and $e_{ij}$'s are errors. Assume, for simplicity, that all the
random effects and errors are i.i.d. $N(0,1)$, so that $\mu$ is the only
unknown parameter. Suppose that $n\rightarrow\infty$, while $m$ is fixed,
say, $m=1$. In this case, $\bar{y}_{1\cdot}=n^{-1}\sum_{j=1}^{n}y_{1j}=
\mu+u_{1}+\bar{v}_{\cdot}+\bar{e}_{1\cdot}$ is a consistent estimator of
the cluster mean, $\mu_{1}=\mu+u_{1}$. On the other hand, the MLE of
$\mu$, which is also $\bar{y}_{1\cdot}$, is inconsistent (because it
converges in probability to $\mu+u_{1}$, which is not equal to $\mu$
with probability one). Note that here the ratio of the number of effects
and number of observations in the cluster is $2/n$. Apparently, this is
sufficient for consistently estimating the mixed effect $\mu+u_{1}$, but
not the fixed effect $\mu$. One might suspect that the case $m=1$ is
somewhat extreme, as $\mu$ and $u_{1}$ are ``inseparable''; but it does
not matter. In fact, for any $m\geq1$, as long as it is fixed, the MLE
of $\mu$ is $\bar{y}_{\cdot\cdot}=(mn)^{-1}\sum_{i=1}^{m}\sum_{j=1}^{n}
y_{ij}=\mu+\bar{u}_{\cdot}+\bar{v}_{\cdot}+\bar{e}_{\cdot\cdot}$, which
converges in probability to $\mu+\bar{u}_{\cdot}$ as $n\rightarrow\infty$,
and $\mu+\bar{u}_{\cdot}\neq\mu$ with probability one. Thus, the only way
that the MLE of $\mu$ can be consistent is to have both $m$ and $n$ go to
$\infty$.

The example also helps to explain why it is necessary to consider the
limiting process $m\wedge n\rightarrow\infty$, instead of something else,
in the open problem. The result of Theorem~\ref{th1} shows that $m\wedge
n\rightarrow\infty$ is also sufficient for the consistency of the MLE.
In fact, from the proof of Theorem~\ref{th1} it follows that, for large $m, n$,
we have with probability tending\vadjust{\goodbreak} to one that the conditional probability
that $p_{\mu}(y)\leq p_{\mu+\varepsilon}(y)$ given $y_{[1]}$ is bounded by
$\gamma^{m\wedge n}$ for some constant $0<\gamma<1$. The corresponding
upper bound under Theorem~\ref{th3} is $e^{-\lambda m_{a}}$ for some
constant $\lambda>0$, where $m_{a}$ is the number of independent
vectors in the subset $y_{[1]}$, and a similar result holds under
Theorem~\ref{th4} with the upper bound being $\exp[-\lambda m_{a}\{1+o(1)\}]$.
The assumption of Theorem~\ref{th3}, namely, (A1), makes sure that $m_{*}
=\min_{1\leq a\leq b}m_{a}\rightarrow\infty$ as the sample size
increases; the assumptions of Theorem~\ref{th4}, namely, (A1) and~(B3),
make sure that, in addition, the $o(1)$ in the above vanishes as
$m_{*}\rightarrow\infty$.

Although estimation of the random effects is not an objective of this
paper, in some cases this is of interest. For example, one may consider
estimating the conditional mean of $y_{ij}$ given $u_{i}$ in the open
problem (which may correspond to the conditional probability of successful
mating with the $i$th female in the salamander problem). As mentioned, the
data are not clustered in this case; in other words, all the data are in
the same cluster, so the ratio of the number of effects over the number of
observations is $(1+m+n)/mn=m^{-1}+n^{-1}+(mn)^{-1}$, which goes to zero
as $m\wedge n\rightarrow\infty$. It is easy to show that $\bar
{y}_{i\cdot}
=n^{-1}\sum_{j=1}^{n}y_{ij}$ is a consistent estimator of $\mathrm{
E}_{\mu}(
y_{ij}|u_{i})=\mathrm{ E}\{h(\mu+u_{i}+\eta)\}$, where $h(x)=e^{x}/(1+e^{x})$
and the (unconditional) expectation is with respect to $\eta\sim N(0,1)$,
$1\leq i\leq m$. Similarly, $\bar{y}_{\cdot j}=m^{-1}\sum_{i=1}^{m}y_{ij}$
is a consistent estimator of $\mathrm{ E}_{\mu}(y_{ij}|v_{j})=\mathrm{
E}\{h(\mu
+\xi+v_{j})\}$, where the (unconditional) expectation is with respect to
$\xi\sim N(0,1)$, $1\leq j\leq n$.

% zodis "Acknowledgments" paliekamas pagal autoriu
\section*{Acknowledgements}
The author wishes to thank all the researchers who have had the opportunity,
and interest, to discuss with the author about the open problem over the
past 15 years, especially those who have spent time thinking about a
solution.
The author is
grateful for the constructive comments from an Associate Editor and two
referees that have led to major improvements of the results and presentation.

\begin{supplement}[id=suppA]
\stitle{Supplementary material}
\slink[doi]{10.1214/13-AOS1084SUPP} %[doi,text={...}] - jei reikia
%suskaldyti doi
\sdatatype{.pdf}
\sfilename{aos1084\_supp.pdf}
\sdescription{The supplementary material is available online at
\href{http://anson.ucdavis.edu/\textasciitilde
jiang/glmmmle.suppl.pdf}{http://anson.ucdavis.edu/}\break
\href{http://anson.ucdavis.edu/\textasciitilde jiang/glmmmle.suppl.pdf}{\textasciitilde jiang/glmmmle.suppl.pdf}.}
\end{supplement}

%In particular, the Wald consistency results were obtained based on a
%referee's suggestion.

% imsref loaded by akundreckaite, 2013-01-31 10:46:03
%
% imsref loaded by akundreckaite, 2013-01-31 13:09:28

\printaddresses


\begin{thebibliography}{32}
% BibTex style file: ims.bst, 2013-01-28
% Default style options (sort=0,type=number).
% Used options (sort=1,type=nameyear).

\bibitem[\protect\citeauthoryear{Booth and Hobert}{1999}]{BooHob99}
\begin{barticle}[auto:STB|2013/01/29|08:09:18]
\bauthor{\bsnm{Booth},~\bfnm{J.~G.}\binits{J.~G.}} \AND
  \bauthor{\bsnm{Hobert},~\bfnm{J.~P.}\binits{J.~P.}}
(\byear{1999}).
\btitle{Maximizing generalized linear mixed model likelihoods with an automated
  Monte Carlo EM algorithm}.
\bjournal{J. Roy. Statist. Soc. B}
\bvolume{61}
\bpages{265--285}.
\bptok{imsref}%
\end{barticle}
\endbibitem

\bibitem[\protect\citeauthoryear{Breslow and Clayton}{1993}]{BreCla93}
\begin{barticle}[auto:STB|2013/01/29|08:09:18]
\bauthor{\bsnm{Breslow},~\bfnm{N.~E.}\binits{N.~E.}} \AND
  \bauthor{\bsnm{Clayton},~\bfnm{D.~G.}\binits{D.~G.}}
(\byear{1993}).
\btitle{Approximate inference in generalized linear mixed models}.
\bjournal{J. Amer. Statist. Assoc.}
\bvolume{88}
\bpages{9--25}.
\bptok{imsref}%
\end{barticle}
\endbibitem

\bibitem[\protect\citeauthoryear{Breslow and Lin}{1995}]{BreLin95}
\begin{barticle}[mr]
\bauthor{\bsnm{Breslow},~\bfnm{Norman~E.}\binits{N.~E.}} \AND
  \bauthor{\bsnm{Lin},~\bfnm{Xihong}\binits{X.}}
(\byear{1995}).
\btitle{Bias correction in generalised linear mixed models with a single
  component of dispersion}.
\bjournal{Biometrika}
\bvolume{82}
\bpages{81--91}.
\bid{doi={10.1093/biomet/82.1.81}, issn={0006-3444}, mr={1332840}}
\bptok{imsref}%
\end{barticle}
\endbibitem

\bibitem[\protect\citeauthoryear{Chen}{1995}]{Che95}
\begin{barticle}[mr]
\bauthor{\bsnm{Chen},~\bfnm{Hung}\binits{H.}}
(\byear{1995}).
\btitle{Asymptotically efficient estimation in semiparametric generalized
  linear models}.
\bjournal{Ann. Statist.}
\bvolume{23}
\bpages{1102--1129}.
\bid{doi={10.1214/aos/1176324700}, issn={0090-5364}, mr={1353497}}
\bptok{imsref}%
\end{barticle}
\endbibitem

\bibitem[\protect\citeauthoryear{Cram{\'e}r}{1946}]{Cra46}
\begin{bbook}[mr]
\bauthor{\bsnm{Cram{\'e}r},~\bfnm{Harald}\binits{H.}}
(\byear{1946}).
\btitle{Mathematical {M}ethods of {S}tatistics}.
\bseries{Princeton Mathematical Series}
\bvolume{9}.
\bpublisher{Princeton Univ. Press}, \blocation{Princeton, NJ}.
\bid{mr={0016588}}
\bptok{imsref}%
\end{bbook}
\endbibitem

\bibitem[\protect\citeauthoryear{Drum and McCullagh}{1993}]{DruMcC93}
\begin{barticle}[mr]
\bauthor{\bsnm{Drum},~\bfnm{Melinda~L.}\binits{M.~L.}} \AND
  \bauthor{\bsnm{McCullagh},~\bfnm{Peter}\binits{P.}}
(\byear{1993}).
\btitle{R{EML} estimation with exact covariance in the logistic mixed model}.
\bjournal{Biometrics}
\bvolume{49}
\bpages{677--689}.
\bid{doi={10.2307/2532189}, issn={0006-341X}, mr={1243484}}
\bptok{imsref}%
\end{barticle}
\endbibitem

\bibitem[\protect\citeauthoryear{Gilks, Richardson and
  Spiegelhalter}{1996}]{GiRiSp96}
\begin{bbook}[mr]
\beditor{\bsnm{Gilks},~\bfnm{W.~R.}\binits{W.~R.}},
  \beditor{\bsnm{Richardson},~\bfnm{S.}\binits{S.}} \AND
  \beditor{\bsnm{Spiegelhalter},~\bfnm{D.~J.}\binits{D.~J.}}, eds.
(\byear{1996}).
\btitle{Markov Chain {M}onte {C}arlo in Practice}.
\bpublisher{Chapman \& Hall}, \blocation{London}.
\bid{mr={1397966}}
\bptok{imsref}%
\end{bbook}
\endbibitem

\bibitem[\protect\citeauthoryear{Jiang}{1998}]{Jia98}
\begin{barticle}[mr]
\bauthor{\bsnm{Jiang},~\bfnm{Jiming}\binits{J.}}
(\byear{1998}).
\btitle{Consistent estimators in generalized linear mixed models}.
\bjournal{J. Amer. Statist. Assoc.}
\bvolume{93}
\bpages{720--729}.
\bid{doi={10.2307/2670122}, issn={0162-1459}, mr={1631373}}
\bptok{imsref}%
\end{barticle}
\endbibitem

\bibitem[\protect\citeauthoryear{Jiang}{1999}]{Jia99}
\begin{barticle}[mr]
\bauthor{\bsnm{Jiang},~\bfnm{Jiming}\binits{J.}}
(\byear{1999}).
\btitle{Conditional inference about generalized linear mixed models}.
\bjournal{Ann. Statist.}
\bvolume{27}
\bpages{1974--2007}.
\bid{doi={10.1214/aos/1017939247}, issn={0090-5364}, mr={1765625}}
\bptok{imsref}%
\end{barticle}
\endbibitem

\bibitem[\protect\citeauthoryear{Jiang}{2007}]{Jia07}
\begin{bbook}[mr]
\bauthor{\bsnm{Jiang},~\bfnm{Jiming}\binits{J.}}
(\byear{2007}).
\btitle{Linear and Generalized Linear Mixed Models and Their Applications}.
\bpublisher{Springer}, \blocation{New York}.
\bid{mr={2308058}}
\bptok{imsref}%
\end{bbook}
\endbibitem

\bibitem[\protect\citeauthoryear{Jiang}{2010}]{Jia10}
\begin{bbook}[mr]
\bauthor{\bsnm{Jiang},~\bfnm{Jiming}\binits{J.}}
(\byear{2010}).
\btitle{Large Sample Techniques for Statistics}.
\bpublisher{Springer}, \blocation{New York}.
\bid{doi={10.1007/978-1-4419-6827-2}, mr={2675055}}
\bptok{imsref}%
\end{bbook}
\endbibitem

\bibitem[\protect\citeauthoryear{Jiang}{2013}]{supp}
\begin{bmisc}[auto]
\bauthor{\bsnm{Jiang},~\bfnm{Jiming}\binits{J.}}
(\byear{2013}).
\bhowpublished{Supplement to ``The subset argument and consistency of MLE in
  GLMM: Answer to an open problem and beyond.''
  DOI:\doiurl{10.1214/13-AOS1084SUPP}.}
\bptok{imsref}%
\end{bmisc}
\endbibitem

\bibitem[\protect\citeauthoryear{Jiang and Zhang}{2001}]{JiaZha01}
\begin{barticle}[mr]
\bauthor{\bsnm{Jiang},~\bfnm{Jiming}\binits{J.}} \AND
  \bauthor{\bsnm{Zhang},~\bfnm{Weihong}\binits{W.}}
(\byear{2001}).
\btitle{Robust estimation in generalised linear mixed models}.
\bjournal{Biometrika}
\bvolume{88}
\bpages{753--765}.
\bid{doi={10.1093/biomet/88.3.753}, issn={0006-3444}, mr={1859407}}
\bptok{imsref}%
\end{barticle}
\endbibitem

\bibitem[\protect\citeauthoryear{Karim and Zeger}{1992}]{KarZeg92}
\begin{barticle}[pbm]
\bauthor{\bsnm{Karim},~\bfnm{M.~R.}\binits{M.~R.}} \AND
  \bauthor{\bsnm{Zeger},~\bfnm{S.~L.}\binits{S.~L.}}
(\byear{1992}).
\btitle{Generalized linear models with random effects; salamander mating
  revisited}.
\bjournal{Biometrics}
\bvolume{48}
\bpages{631--644}.
\bid{issn={0006-341X}, pmid={1637985}}
\bptok{imsref}%
\end{barticle}
\endbibitem

\bibitem[\protect\citeauthoryear{Kullback and Leibler}{1951}]{KulLei51}
\begin{barticle}[mr]
\bauthor{\bsnm{Kullback},~\bfnm{S.}\binits{S.}} \AND
  \bauthor{\bsnm{Leibler},~\bfnm{R.~A.}\binits{R.~A.}}
(\byear{1951}).
\btitle{On information and sufficiency}.
\bjournal{Ann. Math. Statistics}
\bvolume{22}
\bpages{79--86}.
\bid{issn={0003-4851}, mr={0039968}}
\bptok{imsref}%
\end{barticle}
\endbibitem

\bibitem[\protect\citeauthoryear{Lehmann and Casella}{1998}]{LehCas98}
\begin{bbook}[mr]
\bauthor{\bsnm{Lehmann},~\bfnm{E.~L.}\binits{E.~L.}} \AND
  \bauthor{\bsnm{Casella},~\bfnm{George}\binits{G.}}
(\byear{1998}).
\btitle{Theory of Point Estimation},
\bedition{2nd} ed.
\bpublisher{Springer}, \blocation{New York}.
\bid{mr={1639875}}
\bptok{imsref}%
\end{bbook}
\endbibitem

\bibitem[\protect\citeauthoryear{Lele, Dennis and Lutscher}{2007}]{LelDenLut07}
\begin{barticle}[pbm]
\bauthor{\bsnm{Lele},~\bfnm{Subhash~R.}\binits{S.~R.}},
  \bauthor{\bsnm{Dennis},~\bfnm{Brian}\binits{B.}} \AND
  \bauthor{\bsnm{Lutscher},~\bfnm{Frithjof}\binits{F.}}
(\byear{2007}).
\btitle{Data cloning: Easy maximum likelihood estimation for complex ecological
  models using Bayesian Markov chain Monte Carlo methods}.
\bjournal{Ecol. Lett.}
\bvolume{10}
\bpages{551--563}.
\bid{doi={10.1111/j.1461-0248.2007.01047.x}, issn={1461-0248}, pii={ELE1047},
  pmid={17542934}}
\bptok{imsref}%
\end{barticle}
\endbibitem

\bibitem[\protect\citeauthoryear{Lele, Nadeem and
  Schmuland}{2010}]{LelNadSch10}
\begin{barticle}[mr]
\bauthor{\bsnm{Lele},~\bfnm{Subhash~R.}\binits{S.~R.}},
  \bauthor{\bsnm{Nadeem},~\bfnm{Khurram}\binits{K.}} \AND
  \bauthor{\bsnm{Schmuland},~\bfnm{Byron}\binits{B.}}
(\byear{2010}).
\btitle{Estimability and likelihood inference for generalized linear mixed
  models using data cloning}.
\bjournal{J. Amer. Statist. Assoc.}
\bvolume{105}
\bpages{1617--1625}.
\bid{doi={10.1198/jasa.2010.tm09757}, issn={0162-1459}, mr={2796576}}
\bptok{imsref}%
\end{barticle}
\endbibitem

\bibitem[\protect\citeauthoryear{Lin}{1997}]{Lin97}
\begin{barticle}[mr]
\bauthor{\bsnm{Lin},~\bfnm{Xihong}\binits{X.}}
(\byear{1997}).
\btitle{Variance component testing in generalised linear models with random
  effects}.
\bjournal{Biometrika}
\bvolume{84}
\bpages{309--326}.
\bid{doi={10.1093/biomet/84.2.309}, issn={0006-3444}, mr={1467049}}
\bptok{imsref}%
\end{barticle}
\endbibitem

\bibitem[\protect\citeauthoryear{Lin and Breslow}{1996}]{LinBre96}
\begin{barticle}[mr]
\bauthor{\bsnm{Lin},~\bfnm{Xihong}\binits{X.}} \AND
  \bauthor{\bsnm{Breslow},~\bfnm{Norman~E.}\binits{N.~E.}}
(\byear{1996}).
\btitle{Bias correction in generalized linear mixed models with multiple
  components of dispersion}.
\bjournal{J. Amer. Statist. Assoc.}
\bvolume{91}
\bpages{1007--1016}.
\bid{doi={10.2307/2291720}, issn={0162-1459}, mr={1424603}}
\bptok{imsref}%
\end{barticle}
\endbibitem

\bibitem[\protect\citeauthoryear{McCullagh and Nelder}{1989}]{McCNel83}
\begin{bbook}[mr]
\bauthor{\bsnm{McCullagh},~\bfnm{P.}\binits{P.}} \AND
  \bauthor{\bsnm{Nelder},~\bfnm{J.~A.}\binits{J.~A.}}
(\byear{1989}).
\btitle{Generalized Linear Models},
\bedition{2nd} ed.
\bpublisher{Chapman \& Hall}, \blocation{London}.
\bptok{imsref}%
\end{bbook}
\endbibitem

\bibitem[\protect\citeauthoryear{McCulloch}{1994}]{McC94}
\begin{barticle}[auto:STB|2013/01/29|08:09:18]
\bauthor{\bsnm{McCulloch},~\bfnm{C.~E.}\binits{C.~E.}}
(\byear{1994}).
\btitle{Maximum likelihood variance components estimation for binary data}.
\bjournal{J. Amer. Statist. Assoc.}
\bvolume{89}
\bpages{330--335}.
\bptok{imsref}%
\end{barticle}
\endbibitem

\bibitem[\protect\citeauthoryear{McCulloch}{1997}]{McC97}
\begin{barticle}[mr]
\bauthor{\bsnm{McCulloch},~\bfnm{Charles~E.}\binits{C.~E.}}
(\byear{1997}).
\btitle{Maximum likelihood algorithms for generalized linear mixed models}.
\bjournal{J. Amer. Statist. Assoc.}
\bvolume{92}
\bpages{162--170}.
\bid{doi={10.2307/2291460}, issn={0162-1459}, mr={1436105}}
\bptok{imsref}%
\end{barticle}
\endbibitem

\bibitem[\protect\citeauthoryear{McCulloch, Searle and
  Neuhaus}{2008}]{McCSeaNeu08}
\begin{bbook}[mr]
\bauthor{\bsnm{McCulloch},~\bfnm{Charles~E.}\binits{C.~E.}},
  \bauthor{\bsnm{Searle},~\bfnm{Shayle~R.}\binits{S.~R.}} \AND
  \bauthor{\bsnm{Neuhaus},~\bfnm{John~M.}\binits{J.~M.}}
(\byear{2008}).
\btitle{Generalized, Linear, and Mixed Models},
\bedition{2nd} ed.
\bpublisher{Wiley}, \blocation{Hoboken, NJ}.
\bid{mr={2431553}}
\bptok{imsref}%
\end{bbook}
\endbibitem

\bibitem[\protect\citeauthoryear{Spiegelhalter et~al.}{2004}]{Spietal04}
\begin{bmisc}[auto:STB|2013/01/29|08:09:18]
\bauthor{\bsnm{Spiegelhalter},~\bfnm{D.~J.}\binits{D.~J.}},
  \bauthor{\bsnm{Thomas},~\bfnm{A.}\binits{A.}},
  \bauthor{\bsnm{Best},~\bfnm{N.}\binits{N.}} \AND
  \bauthor{\bsnm{Lunn},~\bfnm{D.}\binits{D.}}
(\byear{2004}).
\bhowpublished{WinBUGS version 1.4 user manual. MRC Biostatistics Unit,
  Institute of Public Health, London.}
\bptok{imsref}%
\end{bmisc}
\endbibitem

\bibitem[\protect\citeauthoryear{Sutradhar and Rao}{2003}]{SutRao03}
\begin{barticle}[mr]
\bauthor{\bsnm{Sutradhar},~\bfnm{Brajendra~C.}\binits{B.~C.}} \AND
  \bauthor{\bsnm{Rao},~\bfnm{R.~Prabhakar}\binits{R.~P.}}
(\byear{2003}).
\btitle{On quasi-likelihood inference in generalized linear mixed models with
  two components of dispersion}.
\bjournal{Canad. J. Statist.}
\bvolume{31}
\bpages{415--435}.
\bid{doi={10.2307/3315854}, issn={0319-5724}, mr={2043151}}
\bptok{imsref}%
\end{barticle}
\endbibitem

\bibitem[\protect\citeauthoryear{Torabi}{2012}]{Tor12}
\begin{barticle}[mr]
\bauthor{\bsnm{Torabi},~\bfnm{Mahmoud}\binits{M.}}
(\byear{2012}).
\btitle{Likelihood inference in generalized linear mixed models with two
  components of dispersion using data cloning}.
\bjournal{Comput. Statist. Data Anal.}
\bvolume{56}
\bpages{4259--4265}.
\bid{doi={10.1016/j.csda.2012.04.008}, issn={0167-9473}, mr={2957869}}
\bptok{imsref}%
\end{barticle}
\endbibitem

\bibitem[\protect\citeauthoryear{Wald}{1949}]{Wal49}
\begin{barticle}[mr]
\bauthor{\bsnm{Wald},~\bfnm{Abraham}\binits{A.}}
(\byear{1949}).
\btitle{Note on the consistency of the maximum likelihood estimate}.
\bjournal{Ann. Math. Statistics}
\bvolume{20}
\bpages{595--601}.
\bid{issn={0003-4851}, mr={0032169}}
\bptok{imsref}%
\end{barticle}
\endbibitem

\bibitem[\protect\citeauthoryear{Wang, Tsai and Qu}{2012}]{WanTsaQu12}
\begin{barticle}[auto:STB|2013/01/29|08:09:18]
\bauthor{\bsnm{Wang},~\bfnm{P.}\binits{P.}},
  \bauthor{\bsnm{Tsai},~\bfnm{G.~F.}\binits{G.~F.}} \AND
  \bauthor{\bsnm{Qu},~\bfnm{A.}\binits{A.}}
(\byear{2012}).
\btitle{Conditional inference functions for mixed-effects models with
  unspecified random-effects distribution}.
\bjournal{J. Amer. Statist. Assoc.}
\bvolume{107}
\bpages{725--736}.
\bptok{imsref}%
\end{barticle}
\endbibitem

\bibitem[\protect\citeauthoryear{Wolfowitz}{1949}]{Wol49}
\begin{barticle}[mr]
\bauthor{\bsnm{Wolfowitz},~\bfnm{J.}\binits{J.}}
(\byear{1949}).
\btitle{On {W}ald's proof of the consistency of the maximum likelihood
  estimate}.
\bjournal{Ann. Math. Statistics}
\bvolume{20}
\bpages{601--602}.
\bid{issn={0003-4851}, mr={0032170}}
\bptok{imsref}%
\end{barticle}
\endbibitem

\bibitem[\protect\citeauthoryear{Wu and Liang}{2004}]{WuLia04}
\begin{barticle}[mr]
\bauthor{\bsnm{Wu},~\bfnm{Hulin}\binits{H.}} \AND
  \bauthor{\bsnm{Liang},~\bfnm{Hua}\binits{H.}}
(\byear{2004}).
\btitle{Backfitting random varying-coefficient models with time-dependent
  smoothing covariates}.
\bjournal{Scand. J. Stat.}
\bvolume{31}
\bpages{3--19}.
\bid{doi={10.1111/j.1467-9469.2004.00369.x}, issn={0303-6898}, mr={2042595}}
\bptok{imsref}%
\end{barticle}
\endbibitem

\end{thebibliography}
\end{document}